\documentclass[12pt]{article}
\usepackage{theorem}
\newtheorem{lemma}{Lemma}
\newtheorem{thm}{Theorem}
\usepackage{latexsym}
\usepackage{amssymb}
\usepackage{amsmath}
\usepackage{epic}
\usepackage{pifont}
\font\smallit=cmti10

\begin{document}
\thispagestyle{empty}

\vskip 30pt
\begin{center}
{\bf   ON THE DEGREE OF REGULARITY OF GENERALIZED\\ VAN DER WAERDEN
TRIPLES} \vskip 15pt {\bf Nikos Frantzikinakis}\\
{\smallit Department of Mathematics,
Pennsylvania State University,
University Park, PA 16802}\\{\tt nikos@math.psu.edu}\\
\vskip 10pt
{\bf Bruce Landman}\\
{\smallit Department of Mathematics,  University of West Georgia,
Carrollton, GA 30118}\\{\tt landman@westga.edu}\\
\vskip 10pt
{\bf Aaron Robertson}\\
{\smallit Department of Mathematics,
Colgate University,
Hamilton, NY 13346}\\
{\tt aaron@math.colgate.edu}

\end{center}
\vskip 20pt
\renewcommand{\arraystretch}{1.5}
\centerline{\bf Abstract:} \footnotesize

\noindent Let $1 \leq a \leq b$ be integers. A triple of the form $(x,ax+d,bx+2d)$, where $x,d$ are positive
integers is called an {\em (a,b)-triple}. The {\em degree of regularity} of the family of all $(a,b)$-triples,
denoted dor($a,b)$, is the maximum integer $r$ such that every $r$-coloring of $\mathbb{N}$ admits a
monochromatic $(a,b)$-triple. We settle, in the affirmative, the
conjecture that dor$(a,b) < \infty$ for all $(a,b) \neq (1,1)$. We also
disprove the conjecture that dor($a,b) \in \{1,2,\infty\}$ for all
$(a,b)$.

\normalsize
\section*{\normalsize 1. Introduction}

B.L. van der Waerden [5] proved that for any positive integers $k$ and $r$, there is a positive integer $w(k,r)$
such that any $r$-coloring of $\{1,2,...,w(k,r)\}$ must admit a monochromatic $k$-term arithmetic progression.
In [3], a generalization of van der Waerden's theorem for 3-term arithmetic progressions was investigated.
Namely, for integers $1 \leq a \leq b$, define an {\em $(a,b)$-triple} to be any 3-term sequence of the form
  $(x,ax+d,bx+2d)$, where $x,d$ are positive integers. Taking $a=b=1$ gives a 3-term arithmetic progression, and by
  van der Waerden's theorem the associated van der Waerden number $w(3,r)$ is finite for all $r$.

Throughout this note, we assume that $a$ and $b$ are integers and that $1 \leq a \leq b$. For $r \geq 1$, denote
by $n=n(a,b;r)$ the least positive integer, if it exists, such that every $r$-coloring of $[1,n]$ admits a
monochromatic $(a,b)$-triple. If no such $n$ exists, we write $n(a,b;r)= \infty$. We say that $(a,b)$ is {\it
regular} if $n(a,b;r) < \infty$ for each $r \in \mathbb{N}$. By van der Waerden's theorem $(1,1)$ is regular. If
$(a,b)$ is not regular, the {\it degree of regularity} of $(a,b)$, denoted dor$(a,b)$, is the largest integer
$r$ such that $(a,b)$ is $r$-regular.

In [3], it is shown that for a wide class of pairs $(a,b) \neq (1,1)$, $(a,b)$ is not regular, i.e., dor$(a,b) <
\infty$, and its authors conjectured that, in fact, $(1,1)$ is the {\em only} regular pair. In Section 2 we
confirm this conjecture.

Also in [3], it was shown that \begin{equation} \mbox{dor($a,b$) = 1 if
and only
 if $b=2a$,} \end{equation}
and upper bounds on dor($a,b$) are given for those pairs which are shown  not to be regular. Further,
those authors speculated that dor($a,b) \in \{1,2,\infty\}$ for all pairs
$(a,b)$. In Section 3 we show this conjecture to be false. We also obtain
upper bounds on dor($a,b$) for all $(a,b) \neq (1,1)$, which improve upon
the results of [3], and provide an alternate proof that (1,1) is the only
regular triple.

\section*{\normalsize 2. The Only Regular Triples are Arithmetic Progressions}

In this section we give a short proof which shows that (1,1)-triples are the only regular $(a,b)$-triples. The
proof makes use of Rado's regularity theorem (see [4]) which states, in
particular, that the linear equation
$a_{1}x_{1}+a_{2}x_2+ \cdots + a_{k}x_{k}=0$ has a monochromatic
solution in $\mathbb{N}$  under any finite coloring
of $\mathbb{N}$ if and
only if some nonempty subset of the nonzero coefficients sums to zero. It
also uses the following lemma.

\begin{lemma}
 For all $1\leq a \leq b$, and all $i \geq 1$,
\[ n(a,b;r) \leq n(a+i,b+2i;r),\]
 and hence dor$(a,b) \geq dor(a+i,b+2i)$.
\end{lemma}
Proof. Let $a,b,i$ be given. To prove the lemma, it suffices to show that every $(a+i,b+2i)$-triple is also an
$(a,b)$-triple. Let $X = (x,y,z)$ be an
 $(a+i,b+2i)$-triple. So $y=(a+i)x+d$  and $z=(b+2i)x+2d$ for some $d > 0$.
But then $X$ is also an $(a,b)$ triple, since $y=ax+(ix+d)$  and
 $z = bx +2(ix+d)$. \hfill $\Box$

\begin{thm} Let $1 \leq a \leq b$. If $(a,b) \neq (1,1)$, then $(a,b)$ is not regular.
\end{thm}
{\em Proof.} Since  the triple $\{x,ax+d,bx+2d\}$ satisfies the equation $(2a-b)x - 2y +z = 0$, by Rado's
regularity theorem an $(a,b)$-triple is regular only if $b-2a \in
\{-2,-1,1\}$.  Hence, this leaves three cases to consider: (i) $b=2a+1$,
(ii) $b=2a-1$, and (iii) $b=2a-2$. In [3] it was shown that dor$(1,3)
\leq 3$, dor$(2,3)=2$, and dor$(2,2) \leq 5$. By Lemma 1, these three
facts cover Cases (i), (ii), and (iii), respectively.
  \hfill $\Box$



\noindent {\bf Remark 1} In Section 3 we will show that dor$(2,2) \leq 4$. We see from this fact, the proof of
Theorem 1, and (1), that $2 \leq $ dor$(a,2a-2) \leq 4$ for all $a \geq
2$; that dor($a,2a-1) = 2$ for all $a
\geq 2$; and that $2 \leq$ dor$(a,2a+1) \leq 3$ for all $a \geq 1$.

\section*{\normalsize 3. More on the Degree of Regularity}

Using the Fortran program {\tt AB.f}, available from the third
 author's website\footnote{\tt
http://math.colgate.edu/$\sim$aaron/programs.html}, we have found that
$n(2,2;3) = 88$.  This implies
\begin{equation} \mbox {dor$(2,2) \geq 3$,}
\end{equation} which is a counterexample to
 the suggestion made in [3] that
dor$(a,b) \in \{1,2,\infty\}$ for all $(a,b)$.
The program uses
 a well-known backtracking algorithm (see [4], Algorithm 2, page 31)
which checks that all $3$-colorings of $[1,88]$ contain
a monochromatic $(2,2)$-triple.


Although (2) shows the existence of a pair besides (1,1) whose degree of regularity is
 greater than two, we wonder if dor$(a,b)=2$ for ``almost all" $(a,b)$.
In particular, we pose the following
questions.

\begin{quote}

\noindent {\bf Question 1}
Is it true that dor$(a,b) \leq 2$ whenever $b \neq 2a-2$ and $a \geq 2$?

\vskip 5pt

\noindent {\bf Question 2}
For $b \neq 2a$, are there only a finite number of pairs $(a,b)$
such that \\ \hspace*{66pt} dor$(a,b) \neq 2$?

\end{quote}
 While we do not yet have the answers to these questions, we have
been able to improve the upper bounds for
dor$(a,b)$, as established in [3], for many $(a,b)$-triples.
These new bounds are supplied by the next two
theorems. The proofs of both theorems use the following coloring.

\noindent {\bf Notation}  Let $c \geq 3$ be an integer and let $p = 2 -
\frac{2}{c}$. Denote by $\gamma_c$ the $c$-coloring of
$\mathbb{N}$ defined by coloring, for each $k \geq 0$, the interval
$[p^k, p^{k+1})$ with color $k \, (\bmod\, c)$.

\noindent
{\bf Theorem 2}  Let $a,i,c \in \mathbb{Z}$ such that $a \geq 2$ and $c
\geq 5$.  Define
$p = 2 - \frac{2}{c}$ and let $0 \leq i \leq p^c(p^{c-1}-2)$.
If $a \leq \frac{p^c}{c-1}$,
then dor$(a,a+i) \leq c-1$.

\noindent
{\it Proof.}  We use the $c$-coloring $\gamma_c$.
Assume, for a contradiction, that $\{x,ax+d,(a+i)x+2d\}$
is a monochromatic $(a,a+i)$-triple under $\gamma_c$.
Let $x \in [p^k,p^{k+1})$.  Since $p<2$ and $a \geq 2$, we have
that $ax+d \in [p^{k+cj},p^{k+cj+1})$ for some $j \in \mathbb{N}$.
This gives us that $d > p^{k+cj}-ap^{k+1}$, which, in turn,
gives us $(a+i)x+2d > 2p^{k+cj}-ap^{k+1}+ip^k$.  We now show that
this lower bound is more that $p^{k+cj+1}$:  By choice of
$a$ we have $a \leq p^{c-1}(2-p)$ so that
$2 - \frac{a}{p^{cj}} \geq p$ for all $j \in \mathbb{N}$.
This gives us $2p^{k+cj}-ap^{k+1} > p^{k+cj+1}$ which is
sufficient for all $i \geq 0$.

Next, we will show that $(a+i)x+2d < p^{k+c(j+1)}$.
Since $d<ax+d<p^{k+cj+1}$ and $ix<ip^{k+1}$ it suffices
to show that $2p^{k+cj+1} + ip^{k+1} < p^{k+cj+c}$.
We have $i \leq p^c(p^{c-1}-2)$, which implies that
$2+\frac{i}{p^cj} < p^{c-1}$ for all $j \in \mathbb{N}$, which, in
turn, implies the desired bound.

Hence, we have $p^{k+cj+1}<(a+i)x+2d < p^{k+c(j+1)}$.
By the definition of $\gamma_c$, we see that if $x$ and
$ax+d$ are the same color, then $(a+i)x+2d$ must be
a different color under $\gamma_c$, a contradiction.
\hfill $\Box$

\vskip 10pt

\noindent
{\bf Example} By Theorem 2 and (2), dor$(2,2) \in \{3,4\}$.

\noindent
{\bf Theorem 3}  Let $b,c \in \mathbb{N}$ such that
$b \geq 2$ and $c \geq 5$.  Let $p=2-\frac{2}{c}$.
If $b < \frac{2+p^c}{p}$, then dor$(1,b) \leq c-1$.

\noindent
{\it Proof.} The proof is quite similar to that of Theorem 2. Assume, for a
contradiction, that $\{x,x+d,bx+2d\}$ is
monochromatic  under $\gamma_c$.
Let $x \in [p^k,p^{k+1})$ so that
$bx+2d \in [p^{k+cj},p^{k+cj+1})$ (since $b \geq 2>c$) for some
$j \in \mathbb{N}$.  This gives
$d \geq \frac12 p^{k+cj} - \frac{b}{2}p^{k+1}$ so
that $x+d > p^k +\frac{1}{2}p^{k+cj}-\frac{b}{2}p^{k+1}$.
The condition on $b$ implies that this last bound
is larger than $p^{k+1}$.

We next show that $x+d < p^{k+cj}$.  We have $d < \frac{1}{2}p^{k+cj+1}$
so that $x+d < p^{k+1}+\frac{1}{2}p^{k+cj+1}$.
Since $2<p^{c-1}(2-p)$ for all $c \geq 5$, we have
$p^{k+1}+\frac{1}{2}p^{k+cj+1}<p^{k+cj}$ for all $j \in \mathbb{N}$.
Hence,  $p^{k+1}<x+d<p^{k+cj}$ so that $x+d$ is
not the same color, under $\gamma_c$, as $x$ and $bx+2d$,
a contradiction.
\hfill $\Box$

\noindent {\bf Corollary 1} For $a \geq 1$ and $1 \leq j \leq 5$, dor($a,2a+j) \leq 4$.

\noindent {\it Proof.} This follows from Theorem 3 and Lemma 1.
\hfill $\Box$

 \noindent {\bf Remark 2}  Theorems 2 and 3, along with the following
result from [3], provide an alternate proof of Theorem 1  without the use
of Rado's regularity theorem.

\begin{lemma} Assume $b \geq (2^{3/2}-1)a-2^{3/2}+2$. Then dor$(a,b) \leq \lceil 2 \log_{2} c \rceil,$ where $c
= \lceil b/a \rceil$.
\end{lemma}

Below we give a table showing the known bounds on the degrees of regularity for some small values of $a$ and
$b$. The entries in the table that improve the previously known bounds are marked with *; all others
are from [3]. The improved
bounds for  dor(1,5), dor(1,6),  dor(1,7), dor(1,8), and dor(1,9) follow from Theorem 3;
the upper bound on dor(2,10) follows from Theorem 2; and the
upper bounds on dor(3,4) and dor(3,7) follow from Lemma 1.

$$ \renewcommand{\arraystretch}{.75}
\begin{array}{ll|ll|ll}
(a,b)&
\mathrm{dor}(a,b)&(a,b)&\mathrm{dor}(a,b)&(a,b)&\mathrm{dor}(a,b)
\\
\hline
(1,1)&\infty& (2,2)&3^{*}-4^{*}&(3,3)&2-5 \\
(1,2)&1&(2,3)&2&(3,4)&2-3^{*}\\
(1,3)&2-3& (2,4)& 1 &(3,5)&2  \\
(1,4)&2-4&(2,5)& 2-3 &(3,6)&1\\
(1,5)&2-4^{*}&  (2,6)&2-3&(3,7)&2-3^{*}\\
(1,6)&2-4^{*}& (2,7)&2-4&(3,8)&2-3 \\
(1,7)&2-4^{*}& (2,8)&2-4&(3,9)&2-3\\
(1,8)&2-5^{*}& (2,9)&2-4& (3,10)&2-4\\
(1,9)&2-5^{*}& (2,10)&2-4^{*} & (3,11)& 2-4\\ \hline
\end{array}
$$

\noindent {\bf Acknowledgement} The result that $(1,1)$ is the only regular pair has been independently shown by
Fox and Radoicic [2]. They show that, in fact, dor$(a,b) \leq 23$ for all $(a,b) \neq (1,1)$.

\parskip=5pt
\section*{\normalsize References}
\footnotesize
\parindent=0pt

[1] T. Brown and B. Landman, Monochromatic
Arithmetic Progressions with Large Differences,
{\it Bull. Australian Math. Soc.} {\bf 60} (1999),
21-35.

[2]  J. Fox and R. Radoicic, preprint

[3] B. Landman and A. Robertson,
On Generalized van der Waerden Triples,
{\it Disc. Math.} {\bf 256} (2002), 279-290.

[4] B. Landman and A. Robertson,
Ramsey Theory on the Integers,
STML {\bf 24}, Am. Math. Soc., 2004.

[5] B. L. van der Waerden, Bewis einer Baudetschen Vermutung,
{\it Nieuw. Arch. Wisk.} {\bf 15} (1927), 212-216.

\end{document}